\documentclass[11pt]{article}
\usepackage[centertags]{amsmath}
 \usepackage{amsfonts}
\usepackage{amssymb}
\usepackage{makeidx}
\usepackage{newlfont}
\usepackage{amsthm} 
\usepackage{stmaryrd}
\usepackage[]{titletoc}  
\usepackage{mathrsfs}
\usepackage{stmaryrd}
 \usepackage{graphicx}  
 \usepackage{xypic}
\usepackage{hyperref}

\newlength{\defbaselineskip}
\newcommand{\setlinespacing}[1]%
           {\setlength{\baselineskip}{#1 \defbaselineskip}}

\theoremstyle{plain}
\newtheorem{thm}{Theorem}[section]
\newtheorem{cor}[thm]{Corollary}
\newtheorem{lem}[thm]{Lemma}
\newtheorem{prop}[thm]{Proposition}
\newtheorem{exam}[thm]{Example}
\newtheorem{rem}[thm]{Remark}

\newcommand{\D}{\mathbb{D}}

\newcommand{\bgm}{L_a^2(\mathbb{D}^2)}

\newcommand{\z}{\Omega}
\newcommand{\V}{ \mathcal{V}}
\newcommand{\N}{  \mathbb{Z}_+}
\newcommand{\E}{\mathcal{E}}
\newcommand{\ths}{ z^{a+nk}w^{b+ml} }
\newcommand{\rds}{ reducing subspace}
\newcommand{\ssim}{ \sim_{a,b}}
\newcommand{\nssim}{ \nsim_{a,b}}

\newcommand{\abs}{(a,b)\in \Omega_2}
\newcommand{\abo}{(a,b)\in \Omega_1}

\makeatletter\@addtoreset{equation}{section} \makeatother
\begin{document}
\title {Multiplication operators defined by a class of  polynomials on $L_a^2(\mathbb{D}^2)$}
\author{ Hui Dan\ \ \ \  Hansong Huang}
\date{}
 \maketitle \noindent\textbf{Abstract:} In this paper, we consider those multiplication
 operators $M_p$ on $L_a^2(\mathbb{D}^2)$ defined by a class of  polynomials $p$. Also, this paper  consider
 the reducing subspaces of $M_p$, the von Neumann
 algebra $ \mathcal{W}^*(p)$ generated by $M_p$, and its commutant $\V^*(p)\triangleq \mathcal{W}^*(p)'.$ The structure of
 $\V^*(p)$ is completely determined, along with those reducing subspaces of $M_p.$



\section{Introduction}
 ~~~~It is interesting to consider    reducing subspaces  of
multiplication operators on function spaces and relevant  von
Neumann algebras. This topic  experienced two main phases of
advancement  during the last four decades.
  The first phase  is  mainly concerned with the theme  in the  1970's on
     commutants and reducing subspaces  of   multiplication operators
     on the Hardy space of the unit disk.
 Several remarkable advances  in this period were  made mainly by  Abrahamse and Douglas\cite{AD};
  Cowen \cite{Cow1,Cow2,Cow3};  Baker, Deddens and Ullman\cite{BDU}; Deddens and Wong\cite{DW};
    Nordgren\cite{Nor}; Thomson \cite{T1,T2,T3,T4}, etc.
  It is natural  to consider   the cases on the Bergman space, and the second phase 
 began with  Zhu's conjecture
   on the numbers of minimal reducing subspaces of   multiplication operators  defined by finite Blaschke products\cite{Zhu} in 2000.
    This research  is presently experiencing a period of intense development, during which  
      a lot of  remarkably progress had been made. Several
      notable results mushroomed, concerning a fascinating connection between analysis, geometry, operator and group
theory, see \cite{DSZ, DPW, GH1, GH2,GH3, HSXY, GSZZ, SW,SZZ1,
SZZ2,Zhu}. For the details, one can refer to \cite{Cow4,CW,GH4}  and
\cite{Guo}.

 However,  little is known in the case if the underlying function space is based on a high-dimensional domain, see
\cite{Guo,LZ,SL}.  
Denote by $\mathbb{D}$ the unit disk in the complex plane
$\mathbb{C}$, and $dA(z)$  the normalized area measure over  $\mathbb{D}$.   
Let $\bgm$ denote the Bergman space consisting of all holomorphic
functions over $\D^2 $ which are square integrable with respect to
the normalized volume
measure $dA(z)dA(w)$. 
For any bounded holomorphic function $\phi$ over $\mathbb{D}^2$,
let $M_\phi$ be the multiplication operator defined on the Bergman space $L_a^2(\mathbb{D}^2)$.
 As done in \cite{GH1}-\cite{GH3},    $\mathcal{W}^*(\phi)$ denotes the von Neumann
 algebra generated by $M_\phi$ and  $ \mathcal{V}^*(\phi)\triangleq \mathcal{W}^*(\phi)'$,
 the commutant algebra of $\mathcal{W}^*(\phi)$.
  It is well-known that  $ \mathcal{V}^*(\phi)$  equals the von Neumann algebra
 generated by the orthogonal projections onto $M$, where $M$ run  over all reducing subspaces of $M_\phi$.
Recall that
 a closed subspace $M$ of $L_a^2(\mathbb{D}^2)$ is called \emph{a reducing subspace}  for   $ M_\phi$ if
$M$ is invariant for both $M_\phi$ and $M_\phi^*$. If, in addition,
there is no nonzero reducing subspace $N$ satisfying $N \subsetneqq
M$, then  $M$ is called \emph{minimal}. This is equivalent to say
that $P_M$ is a minimal projection in $\mathcal{V}^*(\phi)$. Two
reducing subspaces $M_1$ and $M_2$ are called \emph{unitarily
equivalent}
 if  there exists a unitary operator  $U$ from $M$ onto
$N$ and $U$ commutes with   $M_\phi$\cite{GH1}.
 In this case, we   write $$M_1  \stackrel{U}\cong M_2.$$ One can show that $M_1  \stackrel{U}\cong M_2$ if
 and only if $P_{M_1}$ and $P_{M_2}$ are equivalent in  $\mathcal{V}^*( \phi)$;
  that is, there is an operator $V$ in $\mathcal{V}^*(\phi)$ such that
$$ V^*V=P_{M_1} \quad and \quad VV^*=P_{M_2}.$$
  In this paper, $$\mathrm{Z}(\phi) \triangleq \mathcal{W}^*(\phi)\cap
  \mathcal{V}^*(\phi),$$
  the center of $\mathcal{W}^*(\phi).$

Put $p(z,w)=z^k+w^l$ where $k,l\geq 1,$ and this paper mainly focus on reducing subspaces and the relevant von Neumann algebras of
 the multiplication operators $M_p.$ Recall that
 $\V^*(p)\triangleq
\{M_p, M_p^*\}'.$ The following is our main theorem.
\begin{thm} Put $p(z,w)=z^k+w^l$. Then the von Neumann algebra $\V^*(p)$                \label{t01}
defined on $\bgm$ is of type I.  Furthermore, there is a positive integer $K$ such that
 $\V^*(p)$ is $*$-isomorphic to the direct sum of   $M_{n_j}(\mathbb{C})$$(1\leq j \leq K),$ with each $n_j=1$ or $n_j=2$.
 Precisely,  put $\delta=GCD(k,l)$, and then $\V^*(p)$ is $*$-isomorphic to
 $$\bigoplus_{i=1}^m  M_{2}(\mathbb{C})\bigoplus  \Big( \bigoplus_{i=1}^{m'}  \mathbb{C} \Big),$$
  where  $m = \frac{\delta^2-\delta}{2}$ and $m'= kl-\delta^2+2\delta $.
\end{thm}
The following consequence is direct, which completely characterizes the commutativity of $\V^*(p)$   algebraically.
\begin{cor} If $p(z,w)=z^k+w^l$ with $k,l\geq 1,$  then the center $\mathrm{Z}(p)$ of   $\V^*(p)$ is nontrivial;
and $\V^*(p)$ is abelian if and only if  $GCD(k,l)=1.$ In this case, $\mathrm{Z}(p)=\V^*(p)$.
\end{cor}
Recall that if $B$ is a finite Blaschke product, then the von
Neumann algebra $\V^*(B)$ defined on $L_a^2(\mathbb{D})$ is always
abelian\cite{DPW}, i.e.  $\mathrm{Z}(B)=\V^*(B)$. On the other hand,
it is interesting that  $\dim \V^*(z^k+w^l)=kl+\delta^2<\infty$ if
$k,l\geq 1,$ though $ \dim \V^*(z^k w^l)=\infty .$ Notice that
studying $   \V^*(z^k w^l)$ is closely related to studying those
reducing subspaces for $M_{z^k w^l}$,
 which is firstly considered  in \cite{LZ},  and completely characterized in \cite{SL}, both on the unweighted and weighted Bergman spaces.
 However, the approach is quite different from that in this paper.
\begin{rem} Consider $p_\alpha (z,w)=z^k+\alpha w^l $ where $|\alpha
|=1$. It is not difficult to see that $M_{p_\alpha} $ is unitarily equivalent to $M_p$,  and hence $\V^*(p)$
 is $*$-isomorphic to $\V^*(p_\alpha)$. Therefore, similar results also hold  as Theorem \ref{t01} and Corollary 1.2.
\end{rem}

 Below, an example will be presented.
\begin{exam} Put $q(z,w)=z^k+w^k$. Then   $\V^*(q)$
is abelian if and only if $k=1.$ In the case of $k=1, $ $\V^*(z+w)$
has exactly $2$ minimal reducing subspace.
\end{exam}

It is worthwhile to mention that the proof of Theorem 1.1   heavily depends on the reducing subspaces of $M_p$, which will be mentioned as
follows. To begin with,  we needs some notations.
 Put
$$\z \triangleq\{(a,b)\in \N^2; 0\leq a \leq k-1, 0 \leq b \leq l-1\}
$$
For each $(a,b)\in \z$, set $$L_{a,b} \triangleq \overline{span \{
z^{a+nk}w^{b+ml};n,m\in \N\}}  ,$$ which is a reducing subspace for
$M_p. $ Furthermore, it is easy to see that
$$ \bgm=\bigoplus_{(a,b)\in \z}  L_{a,b}.$$
One may hope that all these spaces $L_{a,b}$ are minimal reducing
subspaces for $M_p$. In some cases, it is the case. However, it is not always true. 

Put $s \equiv s(a)= \frac{a+1}{k}$ and $t \equiv t(b)=\frac{b+1}{l}$. Notice that
$$s,t\in  (0,1]\cap \mathbb{Q}.$$
Now   $\Omega$ will be divided  into two parts; precisely, put
 $$\z_1=\{(a,b)\in \z: s\neq t\}   \quad \mathrm{and} \quad \z_2=\{(a,b)\in \z: s= t\}.$$
 Both $\z_1$ and $\z_2$ are always nonempty. Let $M_{a,b}^+$ denote
\[\{f\in L_{a,b}: \langle f, e_{a+nk}(z)e_{b+ml}(w) \rangle =\langle f, e_{a+mk}(z)e_{b+nl}(w) \rangle,\   \forall (m,n) \in \N^2
 \}.\]

Similarly, put
\[M_{a,b}^-=\{f\in L_{a,b}: \langle f, e_{a+nk}(z)e_{b+ml}(w) \rangle =-\langle f, e_{a+mk}(z)e_{b+nl}(w) \rangle,\   \forall (a,b) \in \N^2
 \}\]
For each $f\in \bgm, $ $[f]$ denotes the reducing subspace generated
by $f$. Later,  we will  see that   if $(a,b)\in \Omega_2$, then
both $M_{a,b}^+$ and $M_{a,b}^-$ are  reducing subspaces; and
precisely,  $M_{a,b}^+=[z^aw^b]$ and
$M_{a,b}^-=[z^{a+k}w^b-z^aw^{b+l}]$. \vskip2mm
 The following result is an essential integrant of Theorem \ref{t01}. It tells us which are ``obvious" minimal reducing subspaces.
 \begin{thm} If  $(a,b)\in \z_1$, then $L_{a,b}$ is a minimal reducing subspace for $M_p.$
 If   $(a,b) \in \Omega_2$, then both             \label{p01}
 $M_{a,b}^+$ and $M_{a,b}^-$ are minimal reducing subspaces for  $M_p$
 whose direct sum is $L_{a,b}$. The Bergman space $\bgm$ is exactly the
finite  direct sum  of all  these  minimal  reducing subspaces.
\end{thm}
For each subset $F$ of $\bgm,$
let $\mathcal{S}(p)F $ denote the closed linear span of all vectors
$$ \prod_{r=1}^n M_p^{i_r}M_p^{*j_r}h,\, h\in F$$
where $n, i_r,j_r\in\N$ satisfying $$\sum_{r=1}^n ( i_r-j_r)=1.$$ Then by Theorems \ref{t01}
and \ref{p01}, one can show the following consequence, also see Corollary 4.9.
 \begin{cor} For each reducing subspace $M$ for  $M_p$, $[M\ominus \mathcal{S}(p)M ]=M.$\end{cor}

This paper is arranged as follows.

Section 2 gives some notations and some computational lemmas that will be useful later.
Section 3 will determine those ``obvious" minimal reducing subspaces, whose direct sum equals the whole space $\bgm.$
Section 4  will determine which minimal reducing subspaces are unitarily equivalent,
 and then the proof of Theorem \ref{t01} will be presented.

\section{Some preliminary lemmas}
 ~~~~In this section, some computational lemmas will be presented, which will be useful later.

In this paper, we always assume that $k,l,n,m,n',m'\in\N,$ and $p$ always denotes the polynomial
$z^k+w^l.$
Put $$T=M_p^*M_p-M_pM_p^*.$$
Notice that $T=(M_{z^k}^*M_{z^k}-M_{z^k}M_{z^k}^*)+(M_{w^l}^*M_{w^l}-M_{w^l}M_{w^l}^*)$. By careful verification,
 we have  \begin{equation} T \ths=\big(\phi(s,n)+\phi(t,m)\big)                \label{star}
\ths,\end{equation}
where   \[ \phi(u,n)= \Big\{
 \begin{array}{cc}
 \frac{1}{(u+n)(u+n+1)} , \quad   \quad   n\geq 1,\\
\frac{u}{u+1} , \quad  \quad  \quad  \quad  \quad \quad  n=0.  \\
\end{array}
\]
Since $T$ is a diagonal operator, for any reducing subspace $M$  of $M_p,$
 $M$ is necessarily reducing for $T$, and hence $M$ is also reducing for all spectrum projections for $T$.
 This is our start-point for determining the reducing subspaces for $M_p$.

Notice that
$\sigma_p(T)=\{\phi(s,n)+\phi(t,m)\mid (a,b)\in\Omega,  (n,m)\in\N^2\}$. Rewrite
$$\sigma_p(T)=\{\lambda_d\mid d\in\N\},$$ and let
$Q_d$ denote the orthogonal projection onto the  space   $$H_d \triangleq \{ x\in\bgm: Tx=\lambda_dx\}.$$
Then $Q_d\in W^*(p)$. It is clear that
for different $d$, $Q_d$ are orthogonal to each other, and
$$\sum_{d\in\N}Q_d=I,$$
where the left side converges in strong operator topology.
We may define an equivalence on $\z\times \N^2.$ Precisely,
define
$$ (a,b,n,m)\sim (a',b' ,n',m')\Longleftrightarrow
\phi(s,n)+\phi(t,m)=\phi(s',n')+\phi(t',m').$$
In particular, for a fixed pair $(a,b)\in\Omega$,   define an equivalence $\sim_{a,b}$ on
$\N^2$ as follows:
\[(n,m)\sim_{a,b}(n',m')\Longleftrightarrow\phi(s,n)+\phi(t,m)=\phi(s,n')+\phi(t,m').\]
Equivalently,  there exists an integer $d$ such
that
$$T z^{a+nk}w^{b+ml}=\lambda_d z^{a+nk}w^{b+ml} \quad and \quad T z^{a+n'k}w^{b+m'l}=\lambda_d z^{a+n'k}w^{b+m'l}.$$

Put $\mathcal{E}=\{(n,m)\in\N^2\mid n,m\geq1\}$, and write
$\E^c=\N^2\setminus \E$.
\begin{lem}Each equivalence $\Delta$ on $\z\times \N^2$ is a finite set. \label{21}  \end{lem}
\begin{proof} Notice that all  $(a,b)$ above  as are contained in a finite set $\Omega$. To show that each equivalence $\Delta$
on $\z\times \N^2$  is a finite set, it suffices to show that for each fixed pair
$(a,b)\in \Omega$, and $(n,m)\in \N^2$
there are only finitely many $(n',m')$ satisfying $$(n,m)\sim_{a,b}(n',m').$$
After one minute thought, we may assume that both $(n,m)$ and $(n',m') $ lie in $\E.$
In this case,  when $(s(a),t(b))$ is fixed, both $\phi(s,n)$ and $\phi(t,m)$
are strictly decreasing in $n$ or $m$. 
Since $(n',m')\ssim(n,m),$
$$\min\{n',m'\}<\min\{s+n',t+m'\}\leq\max\{s+n,t+m\}<\max\{n+1,m+1\}.$$
If $n'\leq m'$, then $1\leq n'\leq \max\{n+1,m+1\}$, and for each $n'$ there is at most one
integer $m'\geq 1$ satisfying $(n',m')\ssim(n,m)$. Similarly, if $n'> m'$, then $1\leq m'\leq \max\{n+1,m+1\}$,
 and for each $m'$ there is at most one
integer $n'\geq 1$ satisfying $(n',m')\ssim(n,m)$. Thus, there are
 only finitely many $(n',m')$ in $\N^2$ satisfying $$(n,m)\ssim(n',m').$$
The proof is complete.
\end{proof}
\begin{rem} Lemma \ref{21} tells us that each eigen   space $H_d$ for $T$ is of finite dimension,
 where $T=M_p^*M_p-M_pM_p^*$. Therefore, for each $f\in \bgm $, $Q_df$ is always a polynomial, where $Q_d$ is the orthogonal projection onto $H_d$.
\end{rem}
\begin{lem}Suppose that $(n,m)\neq(n',m')$, $n+m=n'+m'$ and $(n,m)\ssim(n',m')$ for some $(a,b)\in \Omega$.
If both $(n,m)$ and $(n',m')$ lie in the same set $\E$ or $\E^c,$
\label{lem22}
 then $\abs$ and $(n,m)=(m',n').$  
 \end{lem}
\begin{proof} Suppose that $(n,m)\neq(n',m')$, $n+m=n'+m'$ and \linebreak $(n,m)\ssim(n',m')$ for some $(a,b)\in \Omega$.
 There are two cases under consideration.
\vskip2mm
\noindent \textbf{Case I}. $(n,m),(n',m')\in \E$. Since $(n,m)\ssim (n',m')$,
\begin{equation} \phi(s,n)+\phi(t,m)= \phi(s,n')+\phi(t,m').\label{01}\end{equation}
Put $\lambda=(s+n)+(m+t)>2,$ and define
$$\rho(x)=\frac{1}{x(x+1)}+\frac{1}{(\lambda-x)(\lambda-x+1)}, \ 0<x< \lambda. $$
Rewrite (\ref{01}) as \begin{equation}\rho(s+n)=\rho(s+n') .\label{00} \end{equation}

First, we show that $\rho$ is strictly decreasing on $(0, \frac{\lambda}{2})$ and  strictly increasing on $ (\frac{\lambda}{2},\lambda).$
To see this, notice that
$$\rho(x)=\frac{1}{x}- \frac{1}{x+1}  +\frac{1}{\lambda-x}-\frac{1}{ \lambda-x+1},\ 0<x< \lambda,$$
and then
$$\rho'(x)=\big[\frac{1}{(x+1)^2}-\frac{1}{x^2}  \big] -\big[\frac{1}{(\lambda-x+1)^2}-\frac{1}{(\lambda-x)^2}\big],\ \ 0<x< \lambda.$$
Since $v(t)=\frac{1}{(t+1)^2}-\frac{1}{t^2}$ is strictly increasing on $(0,+\infty),$
it follows that $\rho'(x)=0$ if and only if $x=\lambda-x$, i.e. $x=\frac{\lambda}{2}.$
Therefore, $\rho'<0$ on $(0, \frac{\lambda}{2})$ and   $\rho'>0$ on $ (\frac{\lambda}{2},\lambda);$  and then
 $\rho$ is strictly decreasing on $(0, \frac{\lambda}{2})$ and  strictly increasing on $ (\frac{\lambda}{2},\lambda).$

Then one can see that for any $x,y\in (0,\lambda),$ $ \rho(x)=\rho(y)$ if and only if $x=y$
or $x+y=\lambda.$
 Since $s+n\neq s+n'$, (\ref{00}) implies that $$s+n+s+n'=\lambda=s+n+m+t,$$ forcing $s+n'=t+m.$
Notice that $s,t\in (0,1]$, and hence $s=t$, which immediately gives $n'=m,$ and thus $(n,m)=(m',n'). $
The identity $s=t$ shows that $\abs.$
\vskip2mm
\noindent \textbf{Case II}.  $(n,m),(n',m')\in \E^c$. In this case, we may assume that $n'=m =0$.  Since $(n,m)\ssim (n',m')$
 and $n=m'$,
$$\frac{1}{(s+n)(s+n+1)}+\frac{t}{t+1}=\frac{s}{s+1}+\frac{1}{(t+n)(t+n+1)}.$$
Then
$$\frac{1}{(s+n)(s+n+1)}+\frac{1}{s+1}=\frac{1}{(t+n)(t+n+1)}+\frac{1}{t+1},$$
forcing $s=t$, because both sides are strictly increasing in $s$ or $t$ on $(0,+\infty).$
Therefore, $\abs.$ Also, it is clear that $(n,m)=(m',n'). $
\end{proof}
\begin{cor} Suppose that $(n,m)\neq(n',m')$, $n+m=n'+m'$, and $(a,b)\in\Omega_1$.
If both $(n,m)$ and $(n',m')$ lie in the same set $\E $ or $\E^c$, then \label{23}
$(n,m)\nssim(n',m').$
\end{cor}
\begin{lem} For each integer $r\geq 1,$  $(r,0)\ssim(0,r)$ if and only if   \linebreak  $(a,b)\not\in\Omega_1$.
If $n\geq 2$ and  $(n+1,0)\ssim(n,1)$,  then $(n ,0)\nssim(n-1,1)$. \label{24}
Similarly, if $m\geq 2$ and  $(0,m+1)\ssim(1,m)$,  then $(0,m)\nssim(1,m-1)$.
\end{lem}
\begin{proof} It is straightforward to check that $(r,0)\ssim(0,r)$ if and only if     $(a,b)\not\in\Omega_1$
 for $r\geq 1.$

We will show that if  $n\geq 2$ and  $(n+1,0)\ssim(n,1)$,  then $(n ,0)\nssim(n-1,1)$.
Assume conversely that $(n ,0)\nssim(n-1,1)$. Then
 $$\frac{1}{(s+n)(s+n+1)}+\frac{t}{t+1}=\frac{1}{(s+n-1)(s+n)}+\frac{1}{(t+1)(t+2)}.$$
Since $(n+1,0)\ssim(n,1)$, $$\frac{1}{(s+n+1)(s+n+2)}+\frac{t}{t+1}=\frac{1}{(s+n)(s+n+1)}+\frac{1}{(t+1)(t+2)},$$
and then
 $$\frac{1}{(s+n)(s+n+1)}- \frac{1}{(s+n+1)(s+n+2)}= \frac{1}{(s+n-1)(s+n)}-\frac{1}{(s+n)(s+n+1)}.$$
 That is, $$\frac{2}{(s+n)(s+n+1)(s+n+2)}=  \frac{2}{(s+n-1)(s+n)(s+n+1)},$$
 which is impossible since $s>0.$ The remaining part is similar.
\end{proof}

\begin{lem}If $(a,b)\in\Omega_2$, then for all $n,m\geq 1, $   $(n+m,0)\nssim(n,m).$ \label{25}
\end{lem}
\begin{proof} Suppose that $(a,b)\in\z_2$ and $n,m\geq 1$. Assume conversely that $(n+m,0)\ssim(n,m).$ That is,
\begin{equation}\frac 1{(s+n+m)(s+n+m+1)}+\frac s{s+1}=\frac 1{(s+n)(s+n+1)}+\frac 1{(s+m)(s+m+1)}. \label{22}
 \end{equation}
Put \begin{eqnarray*} F(x)&=&G(x)\, \big[\frac{1}{(x+n+m)(x+n+m+1)}+\frac{x}{x+1}\\
&-&\frac 1{(x+n)(x+n+1)}-\frac 1{(x+m)(x+m+1)}\big],\end{eqnarray*}
where $G(x)$ is the least  common multiple of $(x+n+m)(x+n+m+1),$
$x+1, (x+n)(x+n+1)  $ and $(x+m)(x+m+1)$. Observe that $F$ is an
integral polynomial whose leading coefficient is $1.$ By basic
algebra, all rational roots of the equation  $F(x)=0$ are integers.
Since $F(s)=0$ and $s\in (0,1]\cap \mathbb{Q}$, then $s=1.$
 However, $s=1$ does not satisfy the equation (\ref{22}); i.e., $F(s)\equiv  F(1)\neq 0$, which is a contradiction.
The proof is complete.\end{proof}
\begin{cor}If $r\geq2$, then \label{26}
$(r,0)$, $(r-1,1)$, $\cdots, (0,r)$  can not belong to a same equivalence with respect to some $(a,b)\in \Omega.$
\end{cor}
\begin{proof} If $(a,b)\in\Omega_1$, then by Lemma \ref{24}
$(r,0)\nssim(0,r)$. If $(a,b)\in\Omega_2$, then Lemma \ref{25} shows that
$(r,0)\nssim(r-1,1),$ completing the proof.
\end{proof}

\section{Minimal  reducing subspaces}
~~~~This section will determine the ``obvious" minimal reducing subspaces,
which is characterized by Theorem  \ref{p01}, restated as below.
 \begin{thm} If  $(a,b) \in \Omega_1$, then $L_{a,b}$ is a minimal reducing subspace for $M_p.$
 If   $(a,b) \in \Omega_2$, then both             \label{p001}
 $M_{a,b}^+$ and $M_{a,b}^-$ are minimal reducing subspaces for  $M_p$
 whose direct sum is $L_{a,b}$. The Bergman space $\bgm$ is exactly the
 finite direct sum of all  these  minimal  reducing subspaces.
\end{thm}
We need some notations. For fixed $(a,b)\in \Omega$ and $r\geq 1,$ put
$$E^r_{a,b}=span\{z^{a+rk}w^b,z^{a+(r-1)k}w^{b+l},\cdots,z^aw^{b+rl}\},$$
and
$E^0_{a,b}= \mathbb{C}z^aw^b.$
Also, notice that
\[\overline{span\{E^r_{a,b}\mid r\in\N\}}=L_{a,b},\]
which is a \rds of $M_p.$
If there is no confusion, we also rewrite $E_r$ for $E^r_{a,b}.$
 \begin{prop} For $ (a,b)\in\Omega$, $[E_r]=L_{a,b},$  $r\geq 1$.   \label{lab}
\end{prop}
 \begin{proof}First, let us make a claim.\vskip2mm

 \noindent \textbf{Claim:}  $[E_r]=[E_{r+1}],\, r\geq 1.$
\vskip2mm
In fact, since $M_p^*z^{a+(r+1)k}w^b=\frac{s+r}{s+r+1}z^{a+rk}w^b$, then
$$z^{a+rk}w^b\in M_p^*(E_{r+1}).$$ Since
$M_p^*z^{a+rk}w^{b+l}=\frac{s+r-1}{s+r}z^{a+(r-1)k}w^{b+l}+\frac t{t+1}z^{a+rk}w^b$,
then $$z^{a+(r-1)k}w^{b+l} \in M_p^*(E_{r+1}).$$ By induction, we have
 $$z^{a+(r-j)k}w^{b+jl} \in M_p^*(E_{r+1}),\ j=0,1,\cdots r,$$ and hence
$E_r\subseteq  M_p^*(E_{r+1}), $  forcing  $[E_r]\subseteq [E_{r+1}].$

To show the inverse inclusion, it suffices to show that $$[E_{r }]\cap  E_{r+1} \supseteq E_{r+1},$$
which is reduced to prove that $$\dim \big( [E_{r }]\cap  E_{r+1}\big) \geq r+2= \dim  E_{r+1}.$$
To see this, we first show that $TM_p(E_r) \nsubseteq M_p(E_r)$ for $r\geq 1$. Otherwise, there is some positive integer $r$
such that $TM_p(E_r) \subseteq M_p(E_r)$. Then for $j=0,1,\cdots, r,$
$$T\big(p z^{a+(r-j)k}w^{b+jl} \big)\in  M_p(E_r) .$$
This, along with (\ref{star}),
 $$T \ths=\big(\phi(s,n)+\phi(t,m)\big)
\ths,\, n,m\in \N,$$  implies   that
there is some nonzero constant $\lambda_j$ satisfying
$$T\big(p z^{a+(r-j)k}w^{b+jl} \big)= \lambda_j \,(p z^{a+(r-j)k}w^{b+jl}). $$
Thus, $$(r+1,0)\ssim(r,1)\ssim\cdots\ssim(0,r+1),$$
which is a contradiction to Corollary \ref{26}. Therefore,  $TM_p(E_r) \nsubseteq M_p(E_r)$.
Also, notice that 
$$\dim M_p(E_r)= \dim  E_r=r+1.$$
and that $$M_p(E_r)\cup TM_p(E_r) \subseteq [E_{r }]\cap  E_{r+1}.$$
Then we have $$\dim  \big( [E_{r }]\cap  E_{r+1} \big)>\dim M_p(E_r)=r+1, $$
forcing $[E_r]\supseteq[E_{r+1}]. $ Therefore,   $[E_r]=[E_{r+1}],$  completing the proof of the claim.

By the above claim,   $[E_1]\supseteq \bigvee_{k\geq 1 }E_r=L_{a,b}.$
On the other hand, $$E_1\subseteq  L_{a,b},$$ which gives $[E_1]\subseteq  L_{a,b} $. Therefore, $$[E_1]=L_{a,b}. $$
Again by the claim, $[E_r]=L_{a,b}, \, r\geq 1. $
\end{proof}
\vskip2mm
By the  proof of Proposition \ref{lab},
\begin{equation}M_p(E_r) \bigvee TM_p(E_r) =E_{r+1},  r\geq 1. \label{eq} \end{equation}
Recall that for each subset $F$ of $\bgm,$
 $\mathcal{S}(p)F $ denotes the closed linear span of all vectors
$$ \prod_{r=1}^n M_p^{i_r}M_p^{*j_r}h,\, h\in F$$
where $n, i_r,j_r\in\N$ satisfying $$\sum_{r=1}^n ( i_r-j_r)=1.$$

\begin{cor} $\dim \bgm \ominus \mathcal{S}( p)\bgm <\infty.$ \end{cor}
Recall that $ E_0= \mathbb{C}z^aw^b $.
 \begin{cor}For $\abo$, $[E_0] \equiv [z^aw^b]=L_{a,b}$. \label{33}
\end{cor}
 \begin{proof}By Proposition \ref{lab}, it suffices to show that $[E_0]\supseteq E_1$.

 Since $\abo,$
$(1,0)\nssim(0,1)$, which implies that $$TM_p(E_0)\nsubseteq M_p(E_0).$$ Then
$\dim [E_0]\cap E_1\geq\dim M_p(E_0)+1=2$, which gives $[E_0]\supseteq E_1$, as desired.
\end{proof}
 By  (\ref{eq}) and Corollary \ref{33}, one can show the following, which will be needed later.
 \begin{lem} For $(a',b')\in \z_1,$ $ L_{a',b'}\ominus \mathcal{S}( p) L_{a',b'}= \mathbb{C} z^{a'}w^{b'}. $ \label{lemma}
 For $(a,b)\in \z_2,$ $ M_{a,b}^+\ominus \mathcal{S}( p) M_{a,b}^+= \mathbb{C} z^aw^b,$ and
   $ M_{a,b}^-\ominus \mathcal{S}( p) M_{a,b}^-= \mathbb{C} z^aw^b(z^k-w^l).$
 \end{lem}
 \begin{prop} For each $\abo$, $L_{a,b}$ is a minimal  \rds \    for $M_p.$   \label{prop2}
\end{prop}
\begin{proof} Assume that $(a,b)\in\Omega_1$. We must prove that for each $ f\in L_{a,b}$ with $f\neq0$, $[f]=L_{a,b};$  by Corollary
\ref{33}, this is equivalent to show that $$z^aw^b\in [f].$$
Since $f\neq 0 $ and $$\sum_{k\in\N} Q_k=I,$$ there exists $ d\in\N$ such that $Q_df\neq0 .$
By Lemma \ref{21} and Remark 2.2, put $$Q_df=\sum_{i=1}^N\alpha_iz^{a+n_ik}w^{b+m_il},$$
where $\alpha_i\neq0$, $\phi(s,n_i)+\phi(t,m_i)=\lambda_d,\,   1\leq i \leq N,$ and
$$n_1+m_1\geq n_2+m_2\geq\cdots\geq m_N+n_N.$$
\noindent \textbf{Case I.} $N=1$ or $N\geq 2$ and $n_1+m_1>n_2+m_2.$ Since
\begin{eqnarray*}M_p^{*n_1+m_1}Q_df &=& M_p^{*n_1+m_1}(\alpha_1z^{a+n_1k}w^{b+m_1l}) \\
 &=& \alpha_1\binom{n_1+m_1}{n_1}\frac{st}{(s+n_1)(t+m_1)}z^aw^b\neq0,
\end{eqnarray*}
  $z^aw^b\in [f].$

\noindent \textbf{Case II.}  $n_1+m_1=n_2+m_2=y\geq1$.
Notice that by Lemma \ref{24} $$(1,0)\nssim(0,1) \quad and  \quad (1,0)\nssim(0,0),$$  and thus
  $y\geq2$.

  Without loss of generality, assume that $N\geq3.$ In this case,  Corollary \ref{23} implies that $$y>m_3+n_3. $$
  Again by Corollary \ref{23}, among $(n_1,m_1)$ and $(n_2,m_2)$, one lies in $\E$ and the other lies in $\E^c.$
  We may assume that  $ (n_1,m_1)=(y,0),$ $n_2\geq 1$ and $m_2\geq 1.$

   Put 
$$\lambda_{d_1}=\phi(s,y-1)+\phi(t,0),\ \, \lambda_{d_2}=\phi(s,n_2-1)+\phi(t,m_2),$$
 and $\lambda_{d_3}=\phi(s,n_2)+\phi(t,m_2-1)$.
There are two situations  under consideration for Case II.

\textbf{1}) \quad $m_2=1$, $n_2=y-1\geq1$.

 Since   $(n_1,m_1)\ssim (n_2,m_2),$ i.e. $(y,0)\ssim (y-1,1)$,
 then Lemma \ref{24} shows that $(y-1,0)\nssim (y-2,1).$ That is,
$ (n_2,m_2-1)\nssim (n_2-1,m_2)$.  Therefore, $$Q_{d_2}z^{a+ (y-1)  k}w^{b }\equiv Q_{d_2}z^{a+ n_2  k}w^{b+(m_2-1)l}=0,$$
and hence
\begin{eqnarray*}Q_{d_2}M_p^*Q_df
=\alpha_2\frac{s+n_2-1}{s+n_2}z^{a+(n_2-1)k}w^{b+m_2l}+\sum_{i=1}^{K'}\alpha_i'z^{a+n'_ik}w^{b+m'_il},\end{eqnarray*}
where $n'_i+m'_i<y-1$, $i=1,2,\cdots,K'.$
Then   $$M_p^{*y-1}Q_{d_2}M_p^*Q_df=\alpha_2\binom{y-1}{m_2}\frac{st}{(s+n_2)(t+m_2)}z^aw^b\neq0,$$
which shows that $z^aw^b\in [f]$, and hence $[z^aw^b]  \subseteq  [f]$.
\vskip2mm
\textbf{2})  $n_2\geq 1$ and $m_2\geq2$.

If $n_2,m_2\geq2$, then by Corollary \ref{23} $$(n_2-1,m_2)\nssim(n_2,m_2-1);$$ and  if
$n_2=1$ and $m_2\geq2$, then $$(y-1,0)\nssim(n_2-1,m_2).$$
In either case, there is a $\lambda_{d_j}(1\leq j \leq 3)$ that is distinct from the other two among $\{\lambda_{d_1},
 \lambda_{d_2}, \lambda_{d_3}\}$. By similar  computations as in \textbf{1}), there is a nonzero constant $\beta_j$ such that
\[M_p^{*y-1}Q_{d_j}M_p^*Q_df=\beta_jz^aw^b.\]
Precisely, $\beta_1=\alpha_1\frac s{s+y-1}$,
$\beta_2=\alpha_2\binom{y-1}{m_2}\frac{st}{(s+n_2)(t+m_2)}$ and
$$\beta_3=\alpha_2\binom{y-1}{n_2}\frac{st}{(s+n_2)(t+m_2)}.$$ Therefore,
$z^aw^b\in [f]$, forcing  $[z^aw^b]  \subseteq  [f]$. The proof is complete.
\end{proof}

Recall that   $M_{a,b}^+$ denotes the reducing  subspace for $M_p$:
\[\{f\in L_{a,b}: \langle f, e_{a+nk}(z)e_{b+ml}(w) \rangle =\langle f, e_{a+mk}(z)e_{b+nl}(w) \rangle,\   \forall (a,b) \in \N^2
 \}\]
 It is not difficult to check that $$\langle f, e_{a+nk}(z)e_{b+ml}(w) \rangle =\langle f, e_{a+mk}(z)e_{b+nl}(w) \rangle$$
  is equivalent to the identity
  $$\langle f, z^{a+nk}w^{b+ml}  \rangle =\langle f, z^{a+mk} w^{b+nl}  \rangle.$$
  Comparing with Corollary \ref{33}, we have the following.
 \begin{lem}For $(a,b)\in\Omega_2$,  \label{35}
$M_{a,b}^+=[z^aw^b]$.
\end{lem}
 \begin{proof} 
Clearly,   $[z^aw^b]\subseteq  M_{a,b}^+$, and it remains to show that
$ M_{a,b}^+\subseteq  [z^aw^b]$.

For this, observe that $M_{a,b}^+\cap E_1 \subseteq [z^aw^b].$ Since
$$M_p(z^{a+k}w^b+z^aw^{b+l})=(z^{a+2k}w^b+z^aw^{b+2l})+2z^{a+k}w^{b+l}$$
and $(0,2)\ssim (2,0)\nssim(1,1)$, then both $z^{a+2k}w^b+z^aw^{b+2l}$ and $z^{a+k}w^{b+l}$ belong to
$M_{a,b}^+\cap E_2$. Since
$\dim (M_{a,b}^+\cap E_2)=2,$
$$M_{a,b}^+\cap E_2\subseteq [z^aw^b].$$ Below, we will use induction to prove that
$$M_{a,b}^+\cap E_r\subseteq [z^aw^b], \ r\geq 1.$$

 By      induction, we assume that $M_{a,b}^+\cap E_r\subseteq [z^aw^b]$,  and we must show that
  $$M_{a,b}^+\cap E_{r+1}\subseteq [z^aw^b];$$
 that is to show $$ z^{a+n k}w^{b+(r-n+1)l}+ z^{a+(r-n+1) k}w^{b+nl} \in[z^aw^b], \, 0\leq n \leq r+1.$$
Notice that either $n\neq \frac{r}{2}$ or  $n-1\neq \frac{r}{2}$.

If $n\neq \frac{r}{2}$, then put $$h(z,w)= z^{a+n k}w^{b+(r-n)l}+ z^{a+(r-n) k}w^{b+nl} \in[z^aw^b].$$
\begin{eqnarray*}M_ph &=&\big( z^{a+(n+1) k}w^{b+(r-n)l}+ z^{a+(r-n) k}w^{b+(n+1)l} \big)+\\
 &+&\big( z^{a+n k}w^{b+(r-n+1)l}+ z^{a+(r-n+1) k}w^{b+nl} \big)\\ &\equiv & f+g.
\end{eqnarray*}
Since $r-n\neq n,$ by Lemmas \ref{lem22} and \ref{25} $(n+1,r-n)\nssim (n,r-n+1)$, which implies that there are two
  constants $\alpha$ and $\beta$ satisfying $\alpha \neq \beta$ and
  $$TM_p h=\alpha f+\beta g.$$
  Also notice that $f+g=M_ph \in [z^aw^b]$, which shows that both $f$ and $g$ lies in $ [z^aw^b]$.

If $n-1\neq  \frac{r}{2}$, then put $$\widetilde{h}(z,w)= z^{a+{n-1} k}w^{b+(r+1-n)l}+ z^{a+(r+1-n) k}w^{b+(n-1)l} \in[z^aw^b].$$
\begin{eqnarray*}M_p\widetilde{h} &=&\big( z^{a+n k}w^{b+(r+1-n)l}+ z^{a+(r+1-n) k}w^{b+nl}
\big)+\\
 &+&\big( z^{a+(n-1) k}w^{b+(r+2-n)l}+ z^{a+(r+2-n) k}w^{b+(n-1)l} \big)\\ &\equiv & \widetilde{f}+\widetilde{g}.
\end{eqnarray*}
Since $r-n+1\neq n-1,$ by Lemmas \ref{lem22} and \ref{25} $$(n,r-n+1)\nssim (n-1,r-n+2),$$ which implies that there are two
  constants $\alpha'$ and $\beta'$ satisfying $\alpha \neq \beta$ and
  $$TM_p h=\alpha' \widetilde{f}+\beta \widetilde{g}.$$
  Also notice that $\widetilde{f}+\widetilde{g}=M_p\widetilde{h} \in [z^aw^b]$,
  and then both $\widetilde{f}$ and $\widetilde{g}$ lies in $ [z^aw^b]$.
%
%

  In either case, we have
$$ z^{a+n k}w^{b+(r-n+1)l}+ z^{a+(r-n+1) k}w^{b+nl} \in[z^aw^b].$$
Therefore, $M_{a,b}^+\cap E_{r+1}\subseteq [z^aw^b]$. The induction is complete and
$$M_{a,b}^+\cap E_r\subseteq [z^aw^b], \ \ r\geq 0.$$
Notice that $M_{a,b}^+ = \bigvee_{k\geq 0}   M_{a,b}^+\cap E_k$. Therefore,
$M_{a,b}^+=[z^aw^b]$.
%
%
%
%
%
%
\end{proof}
 By  similar discussion as the proof of Lemma
\ref{35},  We get the following result.
 \begin{lem} For $(a,b)\in\Omega_2$, \label{prop4}
$M_{a,b}^-=[z^{a }w^b(z^k-w^l)] .$ 
\end{lem}
The following two propostions shows that for  $(a,b)\in\Omega_2$, both $   M_{a,b}^+$ and
$   M_{a,b}^-$ are minimal reducing subspaces.
 \begin{prop}For each  $(a,b)\in\Omega_2$,
$   M_{a,b}^+$ is a minimal \rds. \label{prop5}
\end{prop}
\begin{proof}For  $(a,b)\in\Omega_2$, we must show that  $\forall f\in M_{a,b}^+$ with $f\neq0$,
$[f]=M_{a,b}^+.$ 
In fact, there is an integer $d$ such that $P_d f\neq 0.$ Write
$$Q_df=\sum_{i=1}^K\alpha_i \big( z^{a+n_ik}w^{b+m_il}+ z^{a+m_ik}w^{b+n_il}\big),$$
where   $\alpha_i\neq0$ and $\phi(s,n_i)+\phi(t,m_i)=\lambda_d.$
Without loss of generality, assume that $n_1+m_1\geq
n_2+m_2\geq\cdots\geq m_K+n_K.$ By Lemmas \ref{lem22} and \ref{25}, we may assume that $$n_1+m_1>
n_2+m_2>\cdots> m_K+n_K,$$ and then
$$M_p^{*(n_1+m_1)}Q_d\sum_{i=2}^K\alpha_i \big( z^{a+n_ik}w^{b+m_il}+ z^{a+m_ik}w^{b+n_il}\big) =0.$$
\begin{eqnarray*}M_p^{*(n_1+m_1)}Q_df &=& M_p^{*(n_1+m_1)}Q_d \alpha_1 \big( z^{a+n_1k}w^{b+m_1l}+ z^{a+m_1k}w^{b+n_1l}\big)\\
&=& \frac{2\alpha_1s^2}{(s+n_1)(s+m_1) }\binom{n_1+m_1}{n_1}z^aw^b\neq0,\end{eqnarray*}
forcing $z^aw^b \in [f]$. Then by Lemma \ref{35}, $M_{a,b}^+ =[f]$.
\end{proof}
Also, we have a similar result as Proposition \ref{prop5}.
 \begin{prop}For each  $(a,b)\in\Omega_2$,
$   M_{a,b}^-$ is a minimal \rds. \label{prop6}
\end{prop}
\begin{proof}For  $(a,b)\in\Omega_2$, we must show that  $\forall f\in M_{a,b}^-$ with $f\neq0$,
$[f]=M_{a,b}^-.$ 
In fact, there is an integer $d$ such that $Q_d f\neq 0;$ Write
$$Q_df=\sum_{i=1}^K\alpha_i \big( z^{a+n_ik}w^{b+m_il}- z^{a+m_ik}w^{b+n_il}\big),$$
where   $\alpha_i\neq0$ and $\phi(s,n_i)+\phi(t,m_i)=\lambda_d.$
Without loss of generality, assume that $n_1+m_1\geq
n_2+m_2\geq\cdots\geq m_K+n_K.$ By Lemmas \ref{lem22} and \ref{25}, $n_1+m_1>
n_2+m_2>\cdots> m_K+n_K$, and then
$$M_p^{*(n_1+m_1-1)}Q_d\sum_{i\geq3} \alpha_i \big( z^{a+n_ik}w^{b+m_il}- z^{a+m_ik}w^{b+n_il}\big) =0.$$
Also notice that if   $n_2+m_2=n_1+m_1-1, $ then
$$M_p^{*(n_1+m_1-1)}Q_d  \big( z^{a+n_2k}w^{b+m_2l}- z^{a+m_2k}w^{b+n_2l}\big)=$$$$
 \frac{ \alpha_1s(s+1)}{(s+n_1)(s+m_1) }(\binom{n_1+m_1-1}{n_2}-\binom{n_1+m_1-1}{m_2})z^aw^b=0 .$$
Therefore, $$M_p^{*(n_1+m_1-1)}Q_d\sum_{i\geq 2} \alpha_i \big( z^{a+n_ik}w^{b+m_il}+ z^{a+m_ik}w^{b+n_il}\big) =0,$$ forcing
\begin{eqnarray*}M_p^{*(n_1+m_1-1)}Q_df &=& M_p^{*(n_1+m_1)}Q_d  \alpha_1\big( z^{a+n_1k}w^{b+m_1l}- z^{a+m_1k}w^{b+n_1l}\big)\\
&=& \frac{ \alpha_1s(s+1)}{(s+n_1)(s+m_1) } c(n_1,m_1)  z^aw^b(z^k-w^l) ,\end{eqnarray*}
where $$c(n_1,m_1) \triangleq\binom{n_1+m_1-1}{n_1}-\binom{n_1+m_1-1}{m_1}=\binom{n_1+m_1-1}{n_1}(1-\frac{n_1}{m_1})\neq 0.$$
 Then by Lemma \ref{prop4}, $M_{a,b}^-=[f]$.
\end{proof}
Combing Propositions \ref{prop2}, \ref{prop5} with \ref{prop6}
shows that Theorem \ref{p001} holds.

One lemma  will be established, which will be needed later.
By the proofs of Propositions \ref{prop5} and \ref{prop6}, we have the following. 
\begin{lem}  For  $(a,b)\in\Omega_2$,  suppose that  $f,g\in Q_d \bgm$ and     \label{trick1}
 $$ f=\sum_{i=1}^K\alpha_i \big( z^{a+n_ik}w^{b+m_il}+ z^{a+m_ik}w^{b+n_il}\big),$$
   $$ g=\sum_{i=1}^K\alpha_i \big( z^{a+n_ik}w^{b+m_il}- z^{a+m_ik}w^{b+n_il}\big),$$
where   $\alpha_i\neq0$ and $n_1+m_1\geq
n_2+m_2\geq\cdots\geq m_K+n_K.$ Put $y=n_1+m_1.$ Then  $ M_p^{*y}f\neq 0 $
and  $ M_p^{*(y-1)}g\neq 0.$
  \end{lem}
\section{Unitarily equivalent reducing subspace}
~~~~In this section, we will determine those unitarily equivalent
minimal reducing subspaces for $M_p,$ which are presented in Section
4.  It is shown that all minimal reducing subspaces $M_{a,b}^+$ or
$M_{a,b}^-$ with $a,b\in \Omega_2$
 are  not unitarily equivalent to any other reducing  subspaces.
   However, the case is more difficult for $L_{a,b}$ where  $a,b\in \Omega_1$.

   Recall that for a fixed operator $ M_p$ and two reducing subspaces $M_1$ and $M_2$, if
   $P_{M_1}$ is equivalent
    to $P_{M_2}$ in the von Neumann algebra $\mathcal{V}^*(p)$, then $M_1$ is
   called unitarily equivalent to $M_2.$ Equivalently, there is a unitary operator $U:M_1\to M_2$ which commutes with $M_p$, see \cite{Guo,GH4}.
 \begin{prop}If both $(a,b)$ and  $(a',b')$ lie in $\Omega_2$ and $(a,b)\neq (a',b')$, then
  $M_{a,b}^+$ and $M_{a',b'}^+$ are not unitarily equivalent; and similarly,
  $M_{a,b}^-$ and $M_{a',b'}^-$ are not unitarily equivalent. \label{prop78} 
\end{prop}
 \begin{proof} Assume that both $(a,b)$ and  $(a',b')$ lie in $\Omega_2$ and $(a,b)\neq (a',b')$. To show that
  $M_{a,b}^-$ and $M_{a',b'}^-$ are not unitarily equivalent, rewrite
  $$M_1=M_{a,b}^- \quad and \quad M_2=M_{a',b'}^- .$$ Assume  conversely
  $M_1$ is unitarily equivalent to $M_2$; that is, $P_{M_1}\sim P_{M_2}.$ Then there is a partial isometry $U$ in $\mathcal{V}^*(p)$ such that

   \begin{equation} \label{eqz}
   \left\{ \begin{aligned}
            U^*U=P_{M_1} \quad &  and \quad UU^*=P_{M_2},  \\
                   UM_p=M_pU \quad &  and \quad U^*M_p=M_pU^*.
                             \end{aligned} \right.
                             \end{equation}
   Then $$U(M_1\ominus  \mathcal{S}(p)M_1)= M_2\ominus  \mathcal{S}(p)M_2,$$
    This, combined with Lemma \ref{lemma}, implies that there is a nonzero constant $c$ such
   that $$U(z^{a+k}w^b-z^aw^{b+l})=c(z^{a'+k}w^{b'}-z^{a'}w^{b'+l}). $$
  Also notice that $UT=TU,$ and then $$ \phi(s,1)+\phi(s,0) = \phi(s',1)+\phi(s',0) .$$
   By some computations, we have $s=s',$ forcing $a=a'$, which is a
   contradiction. Therefore, $M_{a,b}^-$ and $M_{a',b'}^-$ are not unitarily equivalent. By similar discussions as above,
    $M_{a,b}^+$ and $M_{a',b'}^+$ are not unitarily equivalent.   \end{proof}
 \begin{prop}Assume that $(a,b),(a',b')\in\Omega_2$. Then  $M_{a,b}^+$ and
$M_{a',b'}^-$  are not unitarily equivalent. 
\end{prop}
\begin{proof} Assume that both $(a,b)$ and  $(a',b')$ lie in $\Omega_2$.  Rewrite
  $$M_1=M_{a,b}^+ \quad and \quad M_2=M_{a',b'}^- .$$ We will show that
  $M_1$ and $M_2$ are not unitarily equivalent.  Assume  conversely
  $M_1$ is unitarily equivalent to $M_2$.
   Then there is a partial isometry $U$ in $\mathcal{V}^*(p)$  satisfying (\ref{eqz}).
   Put $M=(I+  U)M_1$, and  consider $z^{a+2k}w^b+z^aw^{b+2l}\in M_1$.
   There must be an integer $d\in \N$ such that
 $$z^{a+2k}w^b+z^aw^{b+2l}\in Q_d\bgm,$$   where $\lambda_d=\phi(s,2)+\phi(s,0)$.
Put $h_0=U(z^{a+2k}w^b+z^aw^{b+2l})$, and then
$$Th_0=UT(z^{a+2k}w^b+z^aw^{b+2l})=\lambda_d U(z^{a+2k}w^b+z^aw^{b+2l})=\lambda_d h_0,$$
which shows that
$h_0\in Q_d\bgm.$ Now write
\begin{equation}h_0=\sum_{i=1}^K \alpha_i (z^{a'+n_i k}w^{b'+m_il}-z^{a'+m_i k}w^{b'+n_il}), \label{h0}\end{equation}
 where  $\alpha_i\neq 0$ and $n_i>m_i$ for $1\leq i\leq K$, and
$$\lambda_d= \phi(s,n_i)+\phi(t,m_i), \, i=1,\cdots, K.$$
Notice that by Lemmas \ref{lem22} and \ref{25},
 $$y\triangleq n_1+m_1>  n_2+m_2> \cdots >n_K+m_K.$$

Rewrite $g_0=z^{a+2k}w^b+z^aw^{b+2l}$, and then
 $$M_p^{*3}g_0=0 \quad and \quad  M_p^{*2}g_0\neq0.$$  Since $h_0= U g_0$
 and $U$ commutes with $M_p^*, $
we have $$M_p^{*3}h_0=0 \quad and \quad  M_p^{*2}h_0\neq0.$$
By Lemma \ref{trick1},
$M^{*3}_ph_0=0, $ we get $y \leq 3.$
Then in (\ref{h0}),
the possible choices for $(n_i,m_i)$ are:
$$ (3,0) ,  (2,1) ,  (2,0) ,  (1,0) .$$
Notice that no  two of  $ (3,0) ,  (2,1) ,  (2,0) ,  (1,0)  $ are $\sim_{a',b'}$ unitarily equivalent,
and thus only one term $(n_i,m_i)$  appears in (\ref{h0}). Since
 $M_p^{*2}h_0\neq 0,$ either
 $$h_0=\alpha_1(z^{a'+3 k}w^{b'}-z^{a' }w^{b'+3l}) ,$$
 or $$h_0=\alpha_1 (z^{a'+2 k}w^{b'+l}- z^{a'+k }w^{b'+2l}). $$
 Then by some computations,
 $$Q_d M_pM_p^*(g_0+h_0)=\frac{s+1}{s+2}g_0+\alpha h_0\in M,$$
 where $\alpha=\frac{s'+2}{s'+3} $ or $\alpha=\frac{s'}{s'+1} .$

 Since $s', s\in(0,1],$ $\alpha-\frac{s+1}{s+2} \neq 0.$
 Since
 $$(\alpha-\frac{s+1}{s+2})h_0=Q_d M_pM_p^*(g_0+h_0)- \frac{s+1}{s+2}(g_0+h_0) \in M,$$
 $h_0\in M,$ which is a contradiction to  $M\cap M_1=\{0\}.$
 The proof is complete.
  \end{proof}

 \begin{prop}If $(a,b)\in\Omega_1$ and
$(a',b')\in\Omega_2$, then
$L_{a,b}$ is not unitarily equivalent to
$M_{a',b'}^+$.  \label{prop10}
\end{prop}
 \begin{proof}
 Put $M_1=L_{a,b}$ and $M_2= M_{a',b'}^+$.  Assume  conversely
  $M_1$ is unitarily equivalent to $M_2$.
   Then there is a partial isometry $U$ in $\mathcal{V}^*(p)$  satisfying (\ref{eqz}).
Put $M=(I+U)M_1$, and then $P_{M_1}M=M_1.$ Now pick $h_0\in M_2$ such that
$z^{a }w^b+h_0\in M.$
Since
$z^{a }w^b \in M_1\ominus \mathcal{S}(p)M_1$,
$$h_0\in  M_2\ominus \mathcal{S}(p)M_2.$$
Therefore, by Lemma \ref{lemma}  there is a nonzero constant  $c$ such
that
$$h_0=c z^{a'}w^{b'}.$$
Notice that by Lemma \ref{24}, $(0,1)\nssim (1,0).$ Then there are two distinct integers
$d_1$ and $d_2$ such that
$$z^{a+k}w^b \in Q_{d_1}\bgm \quad and \quad z^aw^{b+l} \in   Q_{d_2}\bgm.$$
Since $(a',b')\in\Omega_2$, there is an integer $d$ satisfying
$$z^{a'+k}w^{b'}+z^{a'}w^{b'+l} \in   Q_{d }\bgm.$$
Without loss of generality,   assume that
$d_1\neq d.$ Also remind that $d_1\neq d_2$, and then
$$Q_{d_1}M_p(z^{a }w^b+h_0)=Q_{d_1} \big[z^{a+k}w^b+z^{a }w^{b+l} + c (z^{a'+k}w^{b'}+z^{a'}w^{b'+l})\big] =z^{a+k}w^b.$$
Therefore,
$$ z^{a}w^b \in M_1\cap M.$$
 However, $M_1\cap M=\{0\}$, which is a contradiction. The proof is complete.
%
\end{proof}
Similarly, we have the following.
 \begin{prop}Suppose that $(a,b)\in\Omega_1$ and $(a',b')\in\Omega_2$.
Then $L_{a,b}$ is not unitarily equivalent to $M_{a',b'}^-$ . \label{prop11}
\end{prop}
\begin{proof}
 Put $M_1=L_{a,b}$ and $M_2= M_{a',b'}^-$.    Assume  conversely that
  $M_1$ is unitarily equivalent to $M_2$.
   Then there is a partial isometry $U$ in $\mathcal{V}^*(p)$  satisfying (\ref{eqz}), and thus
  $$\dim \mathcal{S}(p)[M_1 \ominus \mathcal{S}(p)M_1]=\dim \mathcal{S}(p)[M_2 \ominus \mathcal{S}(p)M_2]$$
  By Lemma  \ref{lemma}, $M_1 \ominus \mathcal{S}(p)M_1=\mathbb{C}z^aw^b$
  and $$M_2 \ominus \mathcal{S}(p)M_2=\mathbb{C}z^aw^b(z^k-w^l).$$
  By  the proof of Lemma \ref{33},  $\dim \mathcal{S}(p)[M_1 \ominus \mathcal{S}(p)M_1]=\dim E_1=2.$
  However, $\dim \mathcal{S}(p)[M_2 \ominus \mathcal{S}(p)M_2]=\dim \mathbb{C}z^aw^b(z^{2k}-w^{2l})=1,$
  which is a contradiction. The proof is complete.
\end{proof}

With the help of Propositions \ref{prop78}-\ref{prop11},   by applying operator-theoretic techniques one  gets the following.
\begin{thm}Suppose that $(a,b)\in\Omega_2$ and $M=M_{a,b}^+$ or $M_{a,b}^-$. Then            \label{eqi1}
for any reducing subspace $N,$ $M$ is unitarily equivalent to $N$ if and only if $M=N.$ \end{thm} 



 \begin{prop}Suppose that $(a,b),(a',b')\in\Omega_1$,
$(a,b)\neq(a',b')$, and $(s,t)\neq(t',s')$, then $L_{a,b}$ and \label{prop12}
$L_{a',b'}$ are not unitarily equivalent.
\end{prop}
\begin{proof}Since $(a,b)\neq(a',b')$, then $(s,t)\neq(s',t')$. Also, by assumption $(s,t)\neq(t',s').$
Without loss of generality,     assume that
$s\notin\{s',t'\}$.

Assume  conversely   $M_1$ is unitarily equivalent to $M_2$.
   Then there is a partial isometry $U$ in $\mathcal{V}^*(p)$  satisfying (\ref{eqz}).
Put $M=(I+U)M_1$. Now write $z^{a+k}w^b+h \in M,$ where $h = Uz^{a+k}w^b.$
Since $z^{a+k}w^b \in \mathcal{S}(p)(M_1\ominus \mathcal{S}(p)M_1)$,
 $$h\in \mathcal{S}(p)(M_2\ominus \mathcal{S}(p)M_2).$$
 Therefore, there are three constant $c_0,c_1$ and $c_2$ such
 that $$h =c_0 z^{a'}w^{b'}+c_1 z^{a'+k}w^b+c_2 z^{a'}w^{b'+l}. $$
Since  $h = Uz^{a+k}w^b $ and  $ z^{a+k}w^b  \in Q_d\bgm$ for some integer $d\in \N$ ,    then $$h \in Q_d\bgm.$$
Also notice that no  two of $(0,0)$, $(1,0)$ and $(0,1)$ are $\sim_{a',b'}$ equivalent,
and hence  only one of $c_0,c_1$ and $c_2$ is nonzero. Since $M_p^*z^{a+k}w^b \neq 0$
and \linebreak  $h = Uz^{a+k}w^b,$ we have  $M_p^* h\neq 0,$ forcing $c_0=0.$
Therefore, either \linebreak
$h =c z^{a'+k}w^b$ or $h=c z^{a'}w^{b'+l}$, where $c$ is a nonzero constant.
Then it is not difficult to verify that
$$Q_dM_pM_p^*(z^{a+k}w^b+h )=\frac{ s}{s+1}z^{a+k}w^b+\alpha h,$$
where $\alpha=\frac{s'}{s'+1}$ or $\frac{t'}{t'+1}$.  Notice that $\alpha \neq \frac{ s}{s+1}  $ because
 $s\notin\{s',t'\}$.
Since $$ (Q_dM_pM_p^*-\frac{ s}{s+1}I)(z^{a+k}w^b+h )=(\alpha -\frac{ s}{s+1} )h\in M,$$
$h\in M\cap M_2,$ which is a contradiction to $ M\cap M_2=\{0\}$. Therefore, $L_{a,b}$ is not unitarily equivalent to
$L_{a',b'}$. The proof is complete.
\end{proof}
 \begin{thm}Suppose that $(a,b),(a',b')\in\Omega_1$,
$(a,b)\neq(a',b')$.
Then   $L_{a,b}$ and \label{prop14}
$L_{a',b'}$ are  unitarily equivalent if and only if $(s,t)=(t',s')$.
\end{thm}
 \begin{proof} First assume that  $(s,t)=(t',s')$. In this case, define a unitary operator $U:L_{a,b}\to L_{a',b'}$
 such  that
 $$Ue_{a+mk}(z)e_{b+nl}(w)=e_{a'+nk}(z)e_{b'+ml}(w),\ m,n\in \N.$$
 Then by straightforward computations, one has
 $$M_pUe_{a+mk}(z)e_{b+nl}(w) =U M_p e_{a+mk}(z)e_{b+nl}(w),\, m,n\in \N, $$
 and thus $M_pU=UM_p,$ as desired. That is, $L_{a,b}$ is unitarily equivalent to $L_{a',b'}$.

The converse direction follows directly from Proposition \ref{prop12}.
\end{proof}
Before continuing, let us make an observation.
Let $GCD(k,l)$ denote the greatest common divisor of $k$ and $l$. If $(s,t)=(t',s')$,
 then  $GCD(k,l)\geq 2.$
 Thus, we have the following consequence.
 \begin{cor}Suppose that $(a,b),(a',b')\in\Omega_1$,
  and $GCD(k,l)=1$.
Then   $L_{a,b}$ and
$L_{a',b'}$ are  unitarily equivalent  if and only if $L_{a,b}=L_{a',b'}$.
\end{cor}


 Now we come to the proof of Theorem \ref{t01}.
\vskip2mm
\noindent \textbf{Proof of Theorem \ref{t01}.}
With the help of Theorem \ref{p01}, one can give the proof of first part  as follows.
By Theorem \ref{p01}, there are finitely minimal reducing subspaces $M_1,\cdots, M_K$ of $M_p$ whose direct sum is the whole space
$\bgm.$ Put
$P_j=P_{M_j},\ \, j=1,\cdots, K,$
the orthogonal projection onto $M_j$.
 Then there is a partition of $\{1,2,\cdots,K\}:$ $\Lambda_1, \cdots,\Lambda_{K'}$, satisfying
\begin{itemize}
\item[(i)] If $m,n\in\Lambda_j $ for some $j$, then $P_m$ is   equivalent to $P_n$ in $\V^*(p)$;
\item[(ii)] If $m$ and $n$ lie in different $\Lambda_j$, then $P_m$ is not equivalent to $P_n$.
\end{itemize}
 Then it is easy to see that the von Neumann algebra $\V^*(p)$ is the direct sum of finitely many homogenous algebra, with
each one $*$-isomorphic to $M_n(\mathbb{C})$ where $n\in \N.$ The remaining is to determine the size of $n$ and the numbers of $M_n(\mathbb{C})$.

 By Theorems \ref{eqi1} and \ref{prop14}, the integer $n$ for $M_n(\mathbb{C})$ equals $1$ or $ 2$. By definition, $\Omega_2=\{(a,b): (s= t)\}$.
Put $\delta=GCD(k,l),$ and  write $k=k'\delta$ and $l=l'\delta$.
 Then $\sharp \Omega_2=\delta ,$ and
 there are  $2\delta$ minimal reducing subspaces of the form $M_{a,b}^+$ or $M_{a,b}^-$, and
  $kl- \delta $ minimal reducing subspaces of the form $L_{a,b}.$
 Among these $L_{a,b},$ by Theorem \ref{prop14} $L_{a,b} $ and $L_{a',b'} $ are unitarily equivalent if and only
  if $(s,t)=(t',s')$ or $(a,b)=(a',b').$ In this case, there are $\frac{\delta^2-\delta}{2}$ unitarily equivalent pairs of
  $(L_{a,b} ,L_{a',b'})$.

 Put $m = \frac{\delta^2-\delta}{2}$ and $$m'= 2\delta+ ( kl-\delta) -(\delta^2-\delta)= kl-\delta^2+2\delta .$$
Then $\V^*(p)$ is $*$-isomorphic to
 $$\bigoplus_{i=1}^m  M_{2}(\mathbb{C})\bigoplus  \bigoplus_{i=1}^{m'}  \mathbb{C} ,$$
 as desired. The proof of Theorem \ref{t01} is complete.
$\hfill \square$
\vskip2mm
By Theorems \ref{t01} and \ref{p01}, we have the following consequence.
\begin{cor}Put $p =z^k+  w^l.$ Each nonzero reducing subspace $M$ for $M_p$
must be the direct sum of some (orthogonal) minimal reducing subspaces, each with the following form:
\begin{itemize}
\item[(i)] $M_{a,b}^+$  or $M_{a,b}^-$, where $(a,b)\in \z_2$;
\item[(ii)] $L_{a,b}^-$, where $(a,b)\in \z_1$, and there is no $ (a',b')\in\Omega_1$
satisfying \linebreak $(a,b)\neq(a',b')$ and $(s,t)=(t',s')$;
\item[(iii)]  the reducing subspace $[\alpha z^aw^b+\beta z^{a'}w^{b'}]$,  where $|\alpha|^2+|\beta|^2=1$,
both  $(a,b)$ and $ (a',b')\in\Omega_1$, satisfying $(a,b)\neq(a',b')$ and $(s,t)=(t',s')$.
\end{itemize}
In particular, we have $[M\ominus \mathcal{S}(p)M ]=M.$
\end{cor}

 However, the case on the Hardy space $H^2(\mathbb{D}^2)$ is different. In fact, it is shown in \cite{Da} that
 $\mathcal{V}^*(z^k+w^l)$ is $*$-isomorphic to $ M_{kl}(\mathbb{C})\oplus M_{kl}(\mathbb{C}) $. This is
 interesting not only because the ``size" of the matric  $M_{kl}(\mathbb{C}) $ can be larger than $2$, but also because
  it is different from the case of  $\mathcal{V}^*(B)$ defined on $H^2(\mathbb{D})$, where $B$ is a finite Blaschke product,
   and then  $\mathcal{V}^*(B)$ is   $*$-isomorphic to   $ M_{n}(\mathbb{C})$ with $n=$ order $B $.
  It is worthwhile to mention that  if   $h$ is holomorphic on   $\overline{\mathbb{D}}$,
   then there is a holomorphic function $g  $  on   $\overline{\mathbb{D}}$ and
    a finite Blaschke product $B$ satisfying $h=g\circ B$ and  $\mathcal{V}^*(h)=\mathcal{V}^*(B) $,
    on the Hardy space $H^2(\mathbb{D}),$ the Bergman space $L_a^2(D)$, as well as on the Dirichlet space \cite{T1,T2}.


\noindent{Hui Dan, Department of Mathematics, East China
University of Science and Technology, Shanghai, 200237, China,
E-mail: hiroyuki.sarada@gmail.com

\noindent{Hansong Huang, Department of Mathematics, East China
University of Science and Technology, Shanghai, 200237, China,
E-mail: hshuang@ecust.edu.cn
\end{document}